\documentclass[11pt,draft]{amsart}
\usepackage{setspace,units}
\usepackage{bbm,latexsym,mathrsfs,amscd}

\newtheorem{theorem}{Theorem}[section]
\newtheorem{lemma}[theorem]{Lemma}

\theoremstyle{definition}
\newtheorem{definition}[theorem]{Definition}

\numberwithin{equation}{section}

\RequirePackage{amsmath}
\RequirePackage{amssymb}
\RequirePackage{amsthm}
 \RequirePackage{epsfig}

\begin{document}

\date{}
\title[Real numbers as infinite decimals]
{Real numbers as infinite decimals --- theory and computation }

\author{Nicolas Fardin and Liangpan Li}

\address{Lgt De Lattre de Tassigny, 85000 La Roche-sur-Yon, France}
\email{Nicolas.Fardin@ac-nantes.fr}

\address{Department of Mathematical Sciences,
Loughborough University, LE11 3TU, UK}
 \email{liliangpan@gmail.com}



%\date{July}

\begin{abstract}

In the 16th century, Simon Stevin    initiated a modern approach to decimal
representation of measuring numbers, marking a transition from the
discrete arithmetic  practised by the Greeks to the arithmetic of the
continuum taken for granted today. However,
how to perform arithmetic directly on infinite decimals remains a long-standing problem, which
has seen the popular degeometrisation of  real numbers
since  the first  constructions were published in around 1872.
Our article is devoted to solving this historical problem.
An issue that Hardy called ``a fatal defect"  is also settled.

%Completing Loo-Keng Hua's  approach to the real number system pioneered in 1962,
%this paper  defines arithmetical operations directly on infinite decimals and thus constructs the field of real numbers.
%Basic knowledge about arithmetical operations on terminating decimals is only required.
\end{abstract}

%These operations form the field of real numbers.

\maketitle

\section{Introduction}

In the 16th century, Simon Stevin  (\cite{Katz,OconnorPS,Oconnor}) initiated a modern approach to decimal
representation of measuring numbers, marking a transition from the
discrete arithmetic  practiced by the Greeks to the arithmetic of the
continuum taken for granted today. However,
how to perform arithmetic directly on infinite decimals remains a long-standing problem
(\cite[p. 97]{Birkhoff}, \cite[p. 8]{Brannan}, \cite[p. 11]{Bridger},
 \cite[p. 10]{Courant}, \cite[p. 123]{Gardiner}, \cite[p. 16]{Klazar},
 \cite[p. 80]{Patton}, \cite[p. 400]{Richman}, \cite[pp. 105--106]{Tao}, \cite[p. 739]{Weiss}, \cite{Wu}), which
has seen the popular degeometrisation (\cite{Gamelin,OconnorPS,Oconnor,Richman96}) of  real numbers since
the first constructions were published independently by   M\'{e}ray, Heine,  Cantor and Dedekind in around 1872.
Since then, plenty of attempts have been made to construct the same algebraic structure from
various perspectives, such as   the least upper bound property, the Archimedean property,  equivalence classes,  axiomatic approaches,
the additive group of integers,
 continued fractions, harmonic or alternating series, and so on (\cite{Burn,Weiss}).
% It seems that one would enjoy walking along a circle,  setting starting point and direction arbitrarily.

%``Few mathematical structures have undergone as many revisions or have been presented in as many guises as the real numbers" (\cite{Faltin}).

It is well known that any element of the real number system can be  identified with an infinite decimal,
so why not define arithmetical operations directly on infinite decimals?
There is a long list of mathematicians including Weierstrass and Stolz who prioritise the decimal system over  other constructions since we all learned at school how to perform
arithmetic on  terminating decimals.
But many decimal approaches   lack details, and most of them
are essentially not so different from the earliest theories (see \cite{Klazar} for a literature review).
   Decimal constructions of the real number system are thus rarely seen in modern Mathematical Analysis  textbooks.

Our article is devoted to solving this historical  problem.
Only basic knowledge about elementary arithmetic   on terminating decimals, or equivalently on integers, is required.
We believe our approach  can be used for teaching purposes in  universities, and even in schools in some respects.

The main idea of this paper is as simple as performing arithmetic on integers, but in a slightly different way.
We are usually told to add and multiply numbers from right to left. Why not do so from left to right?

Let us  consider a slightly strange sum $12\bullet\bullet+45\bullet\bullet$,
where each summand is of four digits and the black dots are not specified. Since the sum of digits  in the hundreds place
is 7 not 9, we get
\[12\bullet\bullet+45\bullet\bullet=5\frac{7}{8}\bullet\bullet\]
whose thousands place value 5 is independent of the values of black dots, and hundreds place value can only be 7 or 8.
Here $\frac{7}{8}$ represents the two possible choices, not the fraction that is equal to 0.875.
In exactly the same way, we have
\[0.12\bullet\bullet+0.45\bullet\bullet=0.5\frac{7}{8}\bullet\bullet.\]
Since none of the black dots in the above two examples are deterministic, we can prolong them freely, even to infinity in the second example.
Such an observation led Hua (\cite{Hua}) to define addition on infinite decimals in 1962.  Richman (\cite{Richman}) got the same idea  in 1999.
To summarise, the principle is to do addition locally from right to left but globally from left to right\footnote{To compare,
the addition  on  $p$-adic numbers $\cdots b_2b_1b_0.b_{-1}b_{-2}\cdots b_{-m}$ and $\cdots c_2c_1c_0.c_{-1}c_{-2}\cdots c_{-n}$ works  both locally and globally from right to left.}.

To the best of the authors' knowledge,  we are not aware of any works which
define multiplication in a similar way, although Wu (\cite{Wu}) believes it is doable.

Ahead of stating Hua's definition and our multiplication proposal, we fix some notations.
Our ambient space is (\cite{Courant,Fichtenholz,Kodaira})
\[\mathbb{R}=\{a_0.a_1a_2a_3\cdots\in\mathbb{Z}\times\mathbb{Z}_{10}^{\mathbb{N}}: a_k<9\ \mbox{for infinitely many}\ k\},\]
where $\mathbb{Z}_{10}$ denotes the set $\{0,1,2,\ldots,9\}$.
Note that we exclude  infinite decimals ending in a string of 9s once and for all.
Many readers may be more familiar with the classical decimal system (\cite{Tao})
\[\pm z=\pm a_0.a_1a_2a_3\cdots (z\in\mathbb{R}, z\geq0),\]
but there is no essential difference between the options as long as we primarily focus on arithmetic on non-negative real numbers,
then make a suitable extension. One main reason for choosing $\mathbb{R}$
is because generally given a real number $x$ (or a Dedekind cut, a Cauchy sequence, a map on $\mathbb{Z}$, and so on), we don't need to know whether it is non-negative or negative in advance, but  can find a unique $a_0\in\mathbb{Z}$ so that $a_0\leq x<a_0+1$,
a unique $a_1\in\mathbb{Z}_{10}$ so that $a_0+\frac{a_1}{10}\leq x<a_0+\frac{a_1+1}{10}$, and continue in this way to get
 $x=a_0.a_1a_2a_3\cdots$. A second reason is that the link between the model $\mathbb{R}$ and the earliest theories of real numbers
 can be well explained (see Section \ref{section 4}).

An element $x=a_0.a_1a_2a_3\cdots$  is said to be  terminating  if there exists a non-negative integer $m$ such that $a_k=0$ for $k>m$.
In this case,  write $x=a_0.a_1a_2\cdots a_m$ for simplicity. As usual, $10^{-m}$ is the same as  $0.00\cdots01$ whose last digit 1 is at the $m$-th decimal place. For any element $x=a_0.a_1a_2a_3\cdots$ and any non-negative integer $k$, denote $\theta_k(x)=a_k$, the $k$-th digit of $x$, and
$x_k=a_0.a_1a_2\cdots a_k$,  the truncation of $x$ up to the $k$-th digit.
Defining addition  on  terminating decimals is rather standard.
For example,
\[(-8).765+5.678=(-8)+0.765+5.678=(-8)+6.443=(-2).443.\]
%Hua's definition of addition is as follows.

\begin{definition}[addition, Hua \cite{Hua}]\label{addition} Let $x$, $y$ be elements of $\mathbb{R}$.\\
Case 1: Suppose there exists a non-negative integer $m$ such that $\theta_k(x_k+y_k)=9$ for $k>m$.
Then define
\[x+y=(x_m+y_m)_m+10^{-m}.\]
Case 2: Suppose there exists a sequence of positive integers $k_1<k_2<k_3<\cdots$ such that
$\theta_{k_i}(x_{k_i}+y_{k_i})\neq9$ for $i\in\mathbb{N}$. Then $x+y$ is defined by setting
\[(x+y)_{k_i-1}=(x_{k_i}+y_{k_i})_{k_i-1}\ \ \ (i\in\mathbb{N}).\]
\end{definition}

Note that $\theta_k(x_k+y_k)=9$ if and only if $\theta_k(x)+ \theta_k(y)=9$ and $(x_m+y_m)_m=x_m+y_m$, and both succinct substitutes
  were indeed used in the original definition   in \cite{Hua}.
 But to look for a reasonable multiplication rule,
  one should not study the analogue $\theta_k(x)\times\theta_k(y)$, and the reason will be
  easily seen later on.
An element $a_0.a_1a_2a_3\cdots$ is said to be non-negative or negative if $a_0\geq0$ or $a_0<0$.
Everyone knows about performing multiplication   on non-negative terminating decimals.
Obviously, it suffices to consider  the special case $x+y\leq1$ as the general case is linked by
$xy=10^{2s}\times(\frac{x}{10^s}\times\frac{y}{10^s})$
via a large enough non-negative integer $s$. Our  multiplication proposal  is as follows.

% we  assume without loss of generality that $x+y\leq1$. Here is our proposal.

\begin{definition}\label{defn12} Let $x$, $y$ be non-negative elements of $\mathbb{R}$ such that $x+y\leq1$.\\
Case 1: Suppose there exists a non-negative integer $m$ such that $\theta_k(x_ky_k)=9$ for $k>m$.
Then define
\[xy=(x_my_m)_m+10^{-m}.\]
Case 2: Suppose there exists a sequence of positive integers $k_1<k_2<k_3<\cdots$ such that
$\theta_{k_i}(x_{k_i}y_{k_i})\neq9$ for $i\in\mathbb{N}$. Then $xy$ is defined by setting
\[(xy)_{k_i-1}=(x_{k_i}y_{k_i})_{k_i-1}\ \ \ (i\in\mathbb{N}).\]
\end{definition}

The analogy between both definitions supports the validity of the proposal, and a more convincing explanation
is as follows.
 Let $x,y$ be non-negative such that $x+y\leq1$, and denote $x=x_k+\epsilon_k$, $y=y_k+\delta_k$. Naturally we expect
\begin{align*}
xy=(x_k+\epsilon_k)y=x_ky+\epsilon_ky=x_ky_k+x_k\delta_k+\epsilon_ky,
\end{align*}
which implies that
\begin{equation}\label{F1}0\leq xy-x_ky_k<10^{-k}\end{equation}
as $x+y\leq1$ and the maximum between $\epsilon_k$ and $\delta_k$ is strictly less than $10^{-k}$.
Consequently, if $\theta_k(x_ky_k)\leq8$, then $(xy)_{k-1}=(x_ky_k)_{k-1}$. So the second case of the definition is feasible.
To illustrate  the first case, we study a special case of $m=2$ and $(x_2y_2)_2=0.15$.
It follows from (\ref{F1}) that
\begin{equation}\label{F2}(x_ky_k)_k\leq xy<(x_ky_k)_k+2\cdot10^{-k}.\end{equation}
Letting $k=2$ in (\ref{F2}) gives
$0.15\leq xy<0.17$. Note also $(x_2y_2)_2\leq(x_3y_3)_3$. So
\[0.15=(x_2y_2)_2\leq(x_3y_3)_3\leq xy<0.17.\]
Considering the assumption $\theta_3(x_3y_3)=9$, we get $(x_3y_3)_3=0.159$ or $(x_3y_3)_3=0.169$.
The second situation actually could not happen, because if it did then
\begin{align*}
0.169&=(x_3y_3)_3\leq(x_4y_4)_4\leq xy<0.17,\\
0.1699&=(x_4y_4)_4\leq(x_5y_5)_5\leq xy<0.17,\\
0.16999&=(x_5y_5)_5\leq(x_6y_6)_6\leq xy<0.17,\\
&\cdots\cdots\cdots,
\end{align*}
which implies\footnote{This is an explanation not a proof, as we exclude infinite decimals ending in a string of 9s. One can add them into
the real axis by suitable identification after the field structure is established. } $0.16999\cdots\leq xy<0.17$. This is absurd, so  $(x_3y_3)_3=0.159$.
Similarly,
\begin{align*}
0.159&=(x_3y_3)_3\leq xy<0.161,\\
0.1599&=(x_4y_4)_4\leq xy<0.1601,\\
0.15999&=(x_5y_5)_5\leq xy<0.16001,\\
&\cdots\cdots\cdots,
\end{align*}
which implies $xy=0.16=(x_2y_2)_2+10^{-2}$.

To summarise, the principle is to do multiplication first locally via elementary arithmetic then globally from left to right.
A general definition will be given in the next section.

Finally, we need to show that the above arithmetical operations, no matter how reasonable they may be, form a field.
For whatever reason, many other decimal approaches have stopped here (\cite[p. 16]{Klazar}).
Actually, the details of proving the field structure of various models have drawn
lots of  negative feedback in the past.
Our method is to first establish
\begin{align}
|(x+y)_k-x_k-y_k|&\leq M_110^{-k},\\
|(xy)_k-x_ky_k|&\leq M_210^{-k},
\end{align}
then argue by contradiction.
Here $M_1$ and $M_2$ are positive integers independent of $k$, and  both bounds follow
from the corresponding arithmetical definitions in a few lines.

We end this introduction with several citations of impressions of decimal approaches, which  can be changed after the real number system
is reestablished in Sections \ref{section 2} and \ref{section 3}.
\begin{itemize}
\item ``it is not obvious how to perform arithmetical operations" (Brannan \cite{Brannan})
\item ``any solution involves more and more complications" (Bridger \cite{Bridger})
  \item ``this is not a light task" (Courant \cite{Courant})
  \item ``simply do not work for infinite decimals" (Gardiner \cite{Gardiner})
  \item ``even more tedious to explain multiplication" (Stolz and Gmeiner \cite{Stolz})
  \item ``despite being the most familiar, is actually more complicated" (Tao \cite{Tao})
  \item ``popular approach by novices but is fraught with technical difficulties" (Weiss \cite{Weiss})
\end{itemize}

%We end this introduction with one simple application.
%Note $0.123\times0.456=0.056088$ whose third decimal value is 6 not 9. So according to Definition \ref{defn12},
%\[0.123\bullet\bullet\times0.456\bullet\bullet=0.05\bullet\bullet\bullet\bullet\bullet\bullet\bullet\bullet.\]
%In much the same way, one can use Definition \ref{multiplication} to deduce
%\[1.23\bullet\bullet\times4.56\bullet\bullet=5.\bullet\bullet\bullet\bullet\bullet\bullet\bullet\bullet.\]
%As a consequence,
%\[1.23\bullet\bullet\times4.56\bullet\bullet=10^2\times(0.123\bullet\bullet\times0.456\bullet\bullet)=10^2\times0.05\bullet\bullet\bullet\bullet\bullet\bullet\bullet\bullet
%=5.\bullet\bullet\bullet\bullet\bullet\bullet\bullet\bullet.\]

%the much-talked-about $\sqrt{2}$, and devote the rest sections to justify addition and multiplication and verify arithmetical laws.

%Now the above two definitions look so similar to each other.

%\begin{remark}
%Readers may ask whether we can use Definition \ref{addition} to determine all digits of $\sqrt{2}+\sqrt{3}$.
%So far the answer is no because it is difficult to explicitly determine those  of $\sqrt{2}$ and $\sqrt{3}$
%although the previous example shows how to yield sufficiently many leading digits of $\sqrt{2}$ (and $\sqrt{3}$).
%But we do know the second case of Definition \ref{addition} applies to
%$\sqrt{2}+\sqrt{3}$ because this sum is not a terminating decimal.
%\end{remark}

\section{Justification of definitions}\label{section 2}

\subsection{Addition}

To be precise, by justification of Definition \ref{addition} we mean that in the first case
it is independent of the choices of $m$, and  in the second one
$(x_{k_i}+y_{k_i})_{k_i-1}=(x_n+y_n)_{k_i-1}$
for all $n>k_i$, and $x+y$ is  an element of $\mathbb{R}$. As we will do similar work for multiplication, this is left as an exercise.
It is easy to check that $x+0=0+x=x$ for all $x\in\mathbb{R}$, so 0 is the unital element of the addition.

\subsection{Subtraction}

Hua also gave the following definition of subtraction in \cite{Hua}.

\begin{definition}[subtraction]\label{subtraction} Let $x$, $y$ be elements of $\mathbb{R}$.\\
Case 1: Suppose there exists a non-negative integer $m$ such that $\theta_k(x_k)=\theta_k(y_k)$ for $k>m$.
Then define
$x-y=x_m-y_m.$\\
Case 2: Suppose there exists a sequence of positive integers $k_1<k_2<k_3<\cdots$ such that
$\theta_{k_i}(x_{k_i})\neq\theta_k(y_{k_i})$ for $i\in\mathbb{N}$. Then $x-y$ is defined by setting
\[(x-y)_{k_i-1}=(x_{k_i}-y_{k_i})_{k_i-1}\ \ \ (i\in\mathbb{N}).\]
\end{definition}

We also leave the justification of the subtraction definition  as an exercise.
It is easy to check that $x+(0-x)=(0-x)+x=0$ for all $x\in\mathbb{R}$, so $0-x$, denoted simply by $-x$, is the additive inverse of $x$.
As usual, $x<y$ means $x-y$ is negative, or  more intuitively,
 $x$ appears before  $y$ in the ``\emph{dictionary}" $\mathbb{R}$ (or $\mathbb{Z}\times\mathbb{Z}_{10}^{\mathbb{N}}$). $x$ is positive means $0<x$, and
$x\leq y$ means $y-x$ is non-negative, or equivalently, $x<y$ or $x=y$.

\subsection{Multiplication}

The general definition of multiplication is as follows.

\begin{definition}[multiplication]\label{multiplication}  Let $x,y$ be   elements of $\mathbb{R}$.\\
 (1) Suppose $x,y$ are non-negative. Fix a non-negative integer $s$ such that $x+y\leq10^{s}$.\\
Case 1: Suppose there exists a non-negative integer $m$  such that $\theta_{k}(x_{k+s}y_{k+s})=9$ for $k>m$. Then define
$xy=(x_{m+s}y_{m+s})_m+10^{-m}.$\\
Case 2: Suppose there exists a sequence of positive integers $k_1<k_2<k_3<\cdots$ such that
$\theta_{k_i}(x_{k_i+s}y_{k_i+s})\neq9$ for $i\in\mathbb{N}$. Then $xy$ is defined by setting
\[(xy)_{k_i-1}=(x_{k_i+s}y_{k_i+s})_{k_i-1} \ \ \ (i\in\mathbb{N}).\]
(2) Suppose $x,y$ are negative. Then define
$xy=(-x)(-y).$\\
(3) Suppose only one of $x$ and $y$ is negative. Then define $xy=-(x(-y))$.
\end{definition}

This section is devoted to  justifying  Definition \ref{multiplication} (1), so its assumptions are followed.

Case 1:
First, we claim that for any $n>m$,
 \begin{equation}\label{F44}(x_{n+s}y_{n+s})_m=(x_{m+s}y_{m+s})_m.\end{equation}
 To verify (\ref{F44}), it suffices to consider $n=m+1$,
 and suppose this is the case\footnote{This is crucial as readers could go to revisit the illustrating example from Definition \ref{defn12}.}.
Then
\begin{align*}
x_{n+s}y_{n+s}&=x_{m+s}y_{m+s}+(x_{n+s}-x_{m+s})y_{m+s}+x_{n+s}(y_{n+s}-y_{m+s})\\
&\leq x_{m+s}y_{m+s}+(x_{n+s}+y_{n+s})\cdot\frac{9}{10^{n+s}}
\leq x_{m+s}y_{m+s}+\frac{9}{10^{n}},\end{align*}
where the last inequality is due to  $x_{n+s}+y_{n+s}\leq10^s$.
Thus
 \begin{equation}\label{FAA}(x_{n+s}y_{n+s})_n\leq (x_{m+s}y_{m+s}+\frac{9}{10^{n}})_n
\leq(x_{m+s}y_{m+s})_m+\frac{9}{10^n}+\frac{9}{10^n}.
 \end{equation}
Considering the assumption  $\theta_n(x_{n+s}y_{n+s})=9$, one gets
\begin{equation}\label{FBB}(x_{n+s}y_{n+s})_n=(x_{n+s}y_{n+s})_m+\frac{9}{10^n}.\end{equation}
Combining (\ref{FAA}) and (\ref{FBB}) yields
 \[0\leq(x_{n+s}y_{n+s})_m-(x_{m+s}y_{m+s})_m\leq\frac{9}{10^n}<\frac{1}{10^m},\]
which proves claim (\ref{F44}).
 Next, we claim that for any  $n>m$,
\begin{equation}\label{F43}(x_{n+s}y_{n+s})_n+10^{-n}=(x_{m+s}y_{m+s})_m+10^{-m}.\end{equation}
To verify (\ref{F43}), it suffices to consider $n=m+1$,  and suppose this is the case.
  Recall $\theta_n(x_{n+s}y_{n+s})=9$, so (\ref{F43}) is equivalent to (\ref{F44}). Therefore,
 the  definition is independent of the choice of $m$. On the other hand, it follows from (\ref{F44}) that
  \begin{equation}(x_{n+s}y_{n+s})_m+10^{-m}=(x_{m+s}y_{m+s})_m+10^{-m},\end{equation}
 so  the  definition is also independent of the choice of $s$.

 Case 2:
 We claim that
\begin{equation}\label{F555}(x_ny_n)_{k_i-1}=(x_{k_i+s}y_{k_i+s})_{k_i-1}\end{equation}
for $n>k_i+s$. Similar to the verification of the previous case, one gets
\begin{equation}\label{F45}x_ny_n=x_{k_i+s}y_{k_i+s}+\gamma_n\end{equation}
for some $\gamma_n$ with $0\leq\gamma_n<\frac{1}{10^{k_i}}$.
Considering the assumption $\theta_{k_i}(x_{k_i+s}y_{k_i+s})\leq8$, one has
\begin{equation}\label{F46}x_{k_i+s}y_{k_i+s}=(x_{k_i+s}y_{k_i+s})_{k_i-1}+\epsilon_n\end{equation}
for some $\epsilon_n$ with $0\leq\epsilon_n\leq\frac{9}{10^{k_i}}$. Combining (\ref{F45}) and (\ref{F46}) yields
\[0\leq x_ny_n-(x_{k_i+s}y_{k_i+s})_{k_i-1}<\frac{1}{10^{k_i-1}},\]
which proves claim (\ref{F555}).
Consequently, $xy$ is defined as an element of $\mathbb{Z}\times\mathbb{Z}_{10}^{\mathbb{N}}$.

Case 2 (continued): In this part we continue to show $xy\in\mathbb{R}$. Note that this issue is not so important
as, even if $xy\not\in\mathbb{R}$,  one can identify it with a terminating decimal.
We assume $s=0$ for simplicity, and leave the general case as an exercise.
Suppose $xy\not\in\mathbb{R}$, say for example $xy=a_0.a_1a_2\cdots a_{j-1}999\cdots$, where $j$ is equal to some $k_i$.
  Then, we fix an $l>j$ so that
 $x_n\leq x_j+10^{-j}-10^{-l}$ and $y_n\leq y_j+10^{-j}-10^{-l}$ for all $n\geq l$,
 and the reason why this is possible will be explained later. Thus considering $\theta_j(x_{j}y_{j})\leq8$ and $x_jy_j$ is of at most $2j$
 decimal places to the right of its integer part, one gets
 \begin{align*}
 x_ny_n&\leq x_jy_j+(x_j+y_j)(10^{-j}-10^{-l})+(10^{-j}-10^{-l})^2\\
 &\leq  \big((xy)_{j-1}+8\cdot10^{-j}+(10^{-j}-10^{-2j})\big)+(10^{-j}-10^{-l})+10^{-2j}\\
 &=a_0.a_1a_2\cdots a_{j-1}+10^{-(j-1)}-10^{-l}.
\end{align*}
On the other hand,  fixing an $n>l+1$ with $\theta_n(x_ny_n)\leq8$, one gets
\[x_ny_n\geq(x_ny_n)_{n-1}=(xy)_{n-1}=a_0.a_1a_2\cdots a_{j-1}+10^{-(j-1)}-10^{-(n-1)},\]
which contradicts  the above upper bound for $x_ny_n$. Therefore, we must have $xy\in\mathbb{R}$.
The existence of $l$ can be seen as follows. One can first pick an $l_1>j$ so that $\theta_{l_1}(x)\leq8$, then note for any $n\geq l_1$,
\[x_n\leq x_{j}+\Big(\sum_{i=j+1}^{l_1-1}\frac{9}{10^i}\Big)+\frac{8}{10^{l_1}}+\Big(\sum_{i=l_1+1}^n\frac{9}{10^i}\Big)
=x_j+10^{-j}-10^{-l_1}-10^{-n}\leq x_j+10^{-j}-10^{-l_1}.\]
 Similarly, pick an $l_2>j$ so that $\theta_{l_1}(y)\leq8$, and finally set $l=\max\{l_1,l_2\}$.

We also use $x\times y$ to denote $xy$. It is easy to check that $x\times1=1\times x=x$ for all $x\in\mathbb{R}$, so 1 is the unital element of the multiplication.

\subsection{Reciprocal}

Stevin's idea (\cite{Gowers}) works ideally on defining  reciprocal operation.
Given a positive  $x$, one can find a unique non-negative integer $a_0$ so that
$x\times a_0\leq1< x\times (a_0+1)$,
a unique  $a_1\in\mathbb{Z}_{10}$ so that
$x\times (a_0+\frac{a_1}{10})\leq1< x\times (a_0+\frac{a_1+1}{10})$,
and continue in this way to derive an element $y=a_0.a_1a_2a_3\cdots$ of $\mathbb{Z}\times\mathbb{Z}_{10}^{\mathbb{N}}$.
We leave the verification of the facts $x\in\mathbb{R}$ and $xy=yx=1$ as an exercise.
Therefore, the element $y$,  usually  denoted by $x^{-1}$, is the multiplicative inverse (or reciprocal) of $x$.
The reciprocal of a negative element $z$ is defined to be $-((-z)^{-1})$.

%In this part we only show $y\in\mathbb{R}$ and leave the proof of $xy=1$ as an exercise.
%Suppose $y$ is of the form $y=b_0.b_1b_2\cdots b_m999\cdots$. Then for all $k\geq m$,
%\[xy_k\leq1<x(y_k+10^{-k})=x(y_m+10^{-m}),\]
%which yields
%\[0<x(y_m+10^{-m})-1\leq x(y_k+10^{-k})-xy_k=x\cdot10^{-k}.
%\]
%This is absurd  as the fixed positive constant $x(y_m+10^{-m})-1$ could not be bounded from above by $x\cdot10^{-k}$
%if we let $k$ be large enough. Therefore, $y\in\mathbb{R}$.

\section{Arithmetical laws}\label{section 3}

This section will conclude the proof that $(\mathbb{R},+,\times)$ is a field.
Our strategy  agrees with Conway's suggestion that one should
construct the positive reals before constructing any negative ones (\cite{Richman2}).

\begin{lemma}\label{lemma1}
If $x\neq y$, then  there exists an $l\in\mathbb{N}$ such that  $|x_k- y_k|\geq 10^{-l}$ for  $k> l$.
\end{lemma}

\begin{lemma}\label{lemma2}
Let $x,y$ be elements of $\mathbb{R}$. Then $|(x+y)_k-x_k-y_k|\leq 4\cdot 10^{-k}$ for all $k$.
\end{lemma}

\begin{lemma}\label{lemma3}
Let $x,y$ be non-negative elements of $\mathbb{R}$. Then $|(xy)_k-x_ky_k|\leq M\cdot 10^{-k}$ for all $k$, where $M$ is a positive integer depending only on $x$ and $y$.
\end{lemma}

A proof of Lemma \ref{lemma1} is as follows. Assume without loss of generality that
\[x=a_0.a_1a_2a_3\cdots<y=b_0.b_1b_2b_3\cdots.\] Take first a non-negative integer $m$ such that
$x_m<y_m$, then a positive integer $l>m$ so that $a_l\leq8$. For $k>l$, we have
\begin{align*}
y_k-x_k&\geq y_m-\Big(x_m+\Big(\sum_{i=m+1}^l\frac{a_i}{10^i}\Big)+10^{-l}\Big)
\geq y_m-\Big(x_m+\sum_{i=m+1}^l\frac{9}{10^i}\Big)\\ &=(y_m-x_m-10^{-m})+10^{-l}\geq 10^{-l},
\end{align*}
which finishes the proof. Tracing the justification of Definition \ref{multiplication} can yield a proof of Lemma \ref{lemma3},
so we omit the details.
Lemma \ref{lemma2} can be dealt with in a similar way.

\textbf{Commutative laws}: $x+y=y+x$, $xy=yx$.

These laws are self-evident.

\textbf{Associative laws}: $(x+y)+z=x+(y+z)$, $(xy)z=x(yz)$.

It follows from Lemma \ref{lemma2} that
\[|((x+y)+z)_k-x_k-y_k-z_k|=|((x+y)+z)_k-(x+y)_k-z_k+(x+y)_k-x_k-y_k|\leq 8\cdot 10^{-k}.\]
Similarly,
\[|(x+(y+z))_k-x_k-y_k-z_k|\leq 8\cdot 10^{-k},\]
so
\[|((x+y)+z)_k-(x+(y+z))_k|\leq 16\cdot 10^{-k}.\]
If $(x+y)+z$ and $x+(y+z)$ are not the same, then there exists an  $l\in\mathbb{N}$ such that
\[|((x+y)+z)_k-(x+(y+z))_k|\geq10^{-l}\]
for $k>l$. Consequently, $10^{-l}\leq 16\cdot 10^{-k}$ for $k>l$, which is absurd if we let $k=l+2$. This proves the associative law for addition.
In much the same way, one can establish the associative law for multiplication between three non-negative elements.
The general case is left as an exercise.

\textbf{Distributive law}: $x(y+z)=xy+xz$.\\
Case 1: Suppose $x, y,z$ are non-negative. One can provide a proof that is  similar to that of the associative law for addition.\\
Case 2: Suppose $y$ and $z$ are of the same sign. Then the law follows from Case 1.\\
Case 3: Suppose  $y$ and $z$ are not of the same sign. We can assume  without loss of generality that $y+z$, $-y$, and $z$ are of the same sign.  According to Case 2,
$x(y+z)+x(-y)=xz,$
which yields $x(y+z)=xy+xz$.

To conclude, $(\mathbb{R},+,\times)$ is a field.

\section{Classical theories}\label{section 4}

From now on we will discuss various issues related to real numbers.
Many introductory analysis books (\cite{Fichtenholz,Hardy,Kodaira,Rudin}) choose Dedekind's approach, some (\cite{Tao}) prefer Cantor's theory,
but most avoid a detailed construction. Note that Cantor's work is essentially the same as those of M\'{e}ray and Heine (\cite{Oconnor}).
We should understand that an object could have various characterizations or disguises, and so do real numbers.
Below we discuss Dedekind and Cantor's theories from Stevin's viewpoint of infinite decimal  expansions.

\subsection{Dedekind's theory}
A Dedekind cut $(A|B)$ is  formed of two subsets $A,B$ of $\mathbb{Q}$ such that (\cite{Richman})  $A\cup B=\mathbb{Q}$, $a<b$ for any $a\in A$ and $b\in B$, and $B$ is of no smallest element\footnote{Many authors  replace this uniqueness condition with $A$ having no greatest element (see e.g. \cite{Enderton,Hardy,Kodaira}). If so, then the identification of $(\{x\in\mathbb{Q}:x<1\}|\{x\in\mathbb{Q}:x\geq1\})$ is $0.999\cdots$, which does not belong to $\mathbb{R}$.}.
Given a Dedekind cut $(A|B)$, there  exists a unique integer $a_0$ such that $a_0\in A$, $a_0+1\in B$.
Similarly, there exists a unique integer $a_1\in\mathbb{Z}_{10}$ such that $a_0+\frac{a_1}{10}\in A$, $a_0+\frac{a_1+1}{10}\in B$.
Continuing in this way yields an element $a_0.a_1a_2a_3\cdots$ of $\mathbb{R}$, which is naturally identified with
the cut $(A|B)$. One can easily show such an identification is a bijection.

The biggest disadvantage of Dedekind's approach may be that this language  is rarely used in
 advanced  courses and research activities. According to Gamelin's viewpoint (\cite{Gamelin}),
``It is not clear that even Dedekind grasped the import of what he had done".
%Richman (\cite{Richman2}) also commented: ``Nobody ever wants to go through all the details of the development of the
%real numbers via Dedekind cuts".

%(see also Sections \ref{cantor theory} and \ref{GLBP}, and \cite{Fichtenholz}).

\subsection{Cantor's theory}\label{cantor theory}

A sequence of rational numbers $\{q_n\}_{n=1}^{\infty}$ is said to be Cauchy if
for any $\epsilon>0$, there exists a natural number $N$ (depending on $\epsilon$) such that $|q_m-q_n|<\epsilon$
for all $m,n>N$.
First, show that the given sequence is bounded. By the pigeonhole principle, we then pick an  $a_0\in\mathbb{Z}$ so that infinitely many elements of the sequence lie in $[a_0,a_0+1)$,   an  $a_1\in\mathbb{Z}_{10}$ so that infinitely many elements of the sequence lie in $[a_0+\frac{a_1}{10},a_0+\frac{a_1+1}{10})$. Continuing in this way
yields an element $a_0.a_1a_2a_3\cdots$ of $\mathbb{Z}\mathbb{\times}\mathbb{Z}_{10}^{\mathbb{N}}$.
In most cases this procedure outputs a unique identification element of $\mathbb{R}$.
 Sometimes\footnote{For example, study the sequence $\{1+\frac{(-1)^n}{10^n}\}_{n=1}^{\infty}$.} it could also provide two ``\emph{different}" elements such as $0.999\cdots$ and $1.000\cdots$ so we  identify them,
but three or more identification  elements could never exist all together.
To verify this claim, one needs to explore what Cauchy sequence really means, but this is not so difficult.
Once again, one can easily show that the above identification is a bijection.

Cantor's approach is the first example of completing  metric spaces in functional analysis.
Although he calls Cauchy sequences real numbers,
his  greatest mathematical contributions were inspired by decimal expansions.
For example, Cantor's idea of proving $\mathbb{R}^2$ being the same size as $\mathbb{R}$, which he called continuum, can be described as follows.
Given any two real numbers $x=a_0.a_1a_2a_3\cdots$ and $y=b_0.b_1b_2b_3\cdots$, define
 \[\Psi(x,y)=0.a_1b_1\Box a_2b_2\Box a_3b_3\Box\cdots,\]
 where the sequence of empty boxes is left to encode the integer parts of $x$ and $y$, and there are plenty of ways to do so.
 Obviously, $\Psi$ is injective, so the size of $\mathbb{R}^2$ is not greater than that of $\mathbb{R}$.
  Clearly, the continuum is not greater than the size of $\mathbb{R}^2$ either,
  hence\footnote{This is due to the Cantor-Schr\"{o}der-Bernstein theorem which
  states that if there exist injections $f:A\rightarrow B$ and  $g:B\rightarrow A$, then there exists a bijection $h:A\rightarrow B$ (\cite{Di,Enderton}).} they must be the same.

\section{Completeness and its application to ``a fatal defect" in Analysis}\label{section 5}

In general, completeness means the real axis has no gaps. There are several equivalent ways to characterize
the completeness of $\mathbb{R}$, depending on whether it is regarded as a metric space or a totally ordered set.
If $\mathbb{R}$ is viewed as a metric space, then Cauchy's criterion for convergence is a completeness property;
if it is treated as a totally ordered set, then the least upper bound property plays the same role.
At the end of this section, we will include a standard proof
of the greatest lower bound property, the dual of the  least upper bound property.

Next, we discuss an issue about the concept of angle that many authors have overlooked.
Without a rigorous definition of the length of arcs, any use of
the sine and cosine functions could be flawed.
A common strategy in many books (\cite{Brannan,Chang,Fichtenholz,Hua,Stoll,Zorich}) is to
give a  definition of the length of smooth curves as an application of the integration theory.
The trouble is that they have already used  $\sin x$ or $\cos x$
 as  basic examples before integration.

 Kodaira (\cite{Kodaira}) and Rudin (\cite{Rudin})  observed this issue, but their solutions
involve the complex-valued exponential function.
Courant (\cite[pp. 44--45]{Courant}) gave a traditional definition of $\pi$,
but his area approximation procedure could not  automatically imply that the length of the curve
concatenating two concentric arcs is the sum of those of two pieces.
Without such a done deal, any use of derived trigonometric formulas  (\cite[pp. 81--85]{Demailly}) could be flawed.
Hardy  called the issue ``a fatal defect"  in his course of Pure Mathematics (\cite[p. 316]{Hardy}).
In the following,  a class of curves is introduced to solve this issue.

Let $[a,b]\subset\mathbb{R}$ be a bounded interval. A curve $F=(f_1,f_2):[a,b]\rightarrow\mathbb{R}^2$
is said to be monotone if its coordinate functions $f_1$ and $f_2$ are monotone (increasing or decreasing).
 Since the ``\emph{shortest}" curve connecting two points is a straight line,
we  expect the ``\emph{length}" of $F$, denoted by $\mathscr{L}(F)$, to be an upper bound for that of $\overline{F(a)F(b)}$. Suppose
 we divide $[a,b]$ into two small pieces $[a,b]=[a,c]\cup[c,b]$. Clearly,
\[\mathscr{L}(F)=\mathscr{L}(F|_{[a,c]})+\mathscr{L}(F|_{[c,b]})\geq \mathscr{L}(\overline{F(a)F(c)})+\mathscr{L}(\overline{F(c)F(b)})\geq\mathscr{L}(\overline{F(a)F(b)}),\]
 where the last inequality is due to the sum of two sides of a triangle is not smaller than the third part.
 Continuing in this procedure, the quantity closest to  $\mathscr{L}(F)$ we can find is
  \begin{equation}
 \label{supremum}\sup\Big\{\sum_{i=1}^n\mathscr{L}(\overline{F(c_{i-1})F(c_i)}):a=c_0<c_1<c_2<\cdots<c_n=b,\ \ n\in\mathbb{N}\Big\},
 \end{equation}
 where $\sup (\cdot)$ denotes the least upper bound for a given bounded-above subset of $\mathbb{R}$. Note \begin{align*}
\sum_{i=1}^n\mathscr{L}(\overline{F(c_{i-1})F(c_i)})&=\sum_{i=1}^n\sqrt{(f_1(c_{i-1})-f_1(c_i))^2+(f_2(c_{i-1})-f_2(c_i))^2}\\
&\leq \sum_{i=1}^n\big(|f_1(c_{i-1})-f_1(c_i)|+|f_2(c_{i-1})-f_2(c_i)|\big)\\
&=|f_1(a)-f_1(b)|+|f_2(a)-f_2(b)|,
\end{align*}
where the inequality is due to $\sqrt{x^2+y^2}\leq|x|+|y|$,
and the last equality owes to the monotonicity of $f_1$ and $f_2$. So the existence of (\ref{supremum}) is  guaranteed by the least upper  bound property.

\begin{definition}\label{defn51} Given a monotone curve $F:[a,b]\rightarrow\mathbb{R}^2$, define
its length $\mathscr{L}(F)$ to be the supremum  (\ref{supremum}).
\end{definition}

More generally, a curve $F:[a,b]\rightarrow\mathbb{R}^2$ is said to be rectifiable (\cite{Courant,Rudin}) if (\ref{supremum})
exists (as an element of $\mathbb{R}$), or equivalently if its coordinate functions are of bounded variation.
Jordan's decomposition theorem (\cite[p. 173]{Di}) states that any function of bounded variation
can be written as the difference between two monotone functions,
hence Definition \ref{defn51} is very close to the ultimate scenario of rectifiable curves.
% the flying curve of a bird from its birth to death is rectifiable!
 %Anyone who appreciates this definition could also like Darboux's integration theory (\cite[pp. 208--218]{Stoll}).

 Since the semicircle
$S: t\mapsto (t,\sqrt{1-t^2})\  (-1\leq t\leq 1)$
is two-piece monotone, its length exists and is denoted by $\pi$.
In exactly the same way, one can define the length of subarcs of $S$, which is regarded  the same as the concept of the angle of the corresponding sectors.
One can easily extend the length concept to all arcs, and show that
the length of the curve
concatenating two concentric arcs is the sum of those of two pieces. So we are  free to use any
trigonometric formulas taught in school.

We include a proof of the greatest lower bound property\footnote{Why not the least upper bound property?
Everything is essentially the same except one more identification procedure needs to be included.} whose dual plays a vital role for many authors (\cite{Abian,Gowers,Klazar,Li}) in constructing the real number system.
Let $A=\{a_0^{(\lambda)}.a_1^{(\lambda)}a_2^{(\lambda)}a_3^{(\lambda)}\cdots\in\mathbb{R}:\lambda\in \Lambda\}$
be a non-empty subset
that is bounded from below. Denote $z^{(\lambda)}=a_0^{(\lambda)}.a_1^{(\lambda)}a_2^{(\lambda)}a_3^{(\lambda)}\cdots$.
Pick first the smallest integer $b_0$ from
$\{a_0^{(\lambda)}: z^{(\lambda)}\in A \}$,
then the smallest integer $b_1$ from
$\{a_1^{(\lambda)}: z^{(\lambda)}\in A, a_0^{(\lambda)}=b_0 \}$, and
continue in this way to get an element $b_0.b_1b_2b_3\cdots$, which belongs to  $\mathbb{R}$ and is the greatest lower bound for $A$.
%We refer the readers to \cite[pp. 29--30]{Brannan}

\section{Computational discussions}

In the previous section, we  discussed the concept of rectifiable curves, but did not
explain how to do practical calculation, which is a standard topic in integration.
Considering the field of real numbers has been reestablished in Sections \ref{section 2}  and \ref{section 3},
we study some relevant computational issues.

\subsection{General principles}

Readers may ask which case of the addition (or multiplication) definition happens more frequently.
To answer this question,  knowledge about
 Cantor's continuum or Lebesgue's measure theory (\cite{Stoll}) is required.
(1) Let $x\in\mathbb{R}$ be  fixed, and let $y\in\mathbb{R}$ be arbitrary.
Since the set  of all terminating decimals is countable,
its complement  must be uncountable because  $\mathbb{R}$ is uncountable.
So in this sense, the second case is more likely to be applied to compute  $x+y$. (2)
Let $x$ and $y$ be arbitrary. If the cardinality criterion is replaced by measure, then the answer is the same;
if not, we first need go back to Cantor's paradise
and understand $\mathbb{R}^2$ as being the same size as $\mathbb{R}$ (see Section \ref{cantor theory}),
then give a reasonable explanation.

\subsection{Optimal bounds}

Hua's definition of addition implies that
\[(x_{k+1}+y_{k+1})_k\leq (x+y)_k\leq(x_{k+1}+y_{k+1})_k+10^{-k},\]
which is tighter than the more frequently used
\[x_k+y_k\leq(x+y)_k\leq x_k+y_k+2\cdot10^{-k}.\]
Thus $(x+y)_k$ has only two choices: one is $(x_{k+1}+y_{k+1})_k$; the other is $(x_{k+1}+y_{k+1})_k+10^{-k}$.
 The chance (or probability) of getting the first one is 0.95, and its proof is left as an exercise.
To be clear, given the integer parts of $x$ and $y$ have nothing to do with the purpose,  the sample space  is taken to be $[0,1)\times[0,1)$.

Next, consider the multiplication and let  $\Omega=\{(x,y)\in\mathbb{R}^2:x\geq0, y\geq0, x+y\leq1\}$
be the sample space with doubled Lebesgue measure. According to Definition \ref{defn12}, $(xy)_k$ has only two choices: one is $(x_{k+1}y_{k+1})_k$; the other is $(x_{k+1}y_{k+1})_k+10^{-k}$. Since the chance of getting the first choice on each vertical segment $\{(x,y)\in\Omega: 0\leq y\leq 1-x\}$
with $x$ being fixed is not less than 0.9, so is the chance of the same choice on $\Omega$  by Fubini's theorem.

%Similarly,
% $(x-y)_k$ has  two choices, one is $(x_{k+1}-y_{k+1})_k-10^{-k}$, the other is $(x_{k+1}-y_{k+1})_k$, and
%the chance of getting the second one is also 0.95.

%According to Section \ref{section 5}, $\pi$ is the length of a semicircle of radius one.
%But to use Hua's definition to calculate the digits of the sum $\pi+\sqrt{2}$,
%those of both $\pi$ and $\sqrt{2}$ have to be given in advance.
%This is actually a tough task because computing millions, trillions, or any prescribed finite number of digits of  both numbers
%is not equivalent to knowing all of their digits. For example,  whether
%

%We should point out that these  principles only provide a general scenario, and for any specified
%sums or products, one has to study the operations carefully, and sometimes it is very difficult to
%know which case could happen.

\subsection{Turing's computable  numbers}

A real number consists of an infinitely long string of elements of $\mathbb{Z}_{10}$, so in principle its digits carry an infinite amount of
information.
According to Turing (\cite{Turing}, see also \cite{Minshy}), a \textsf{real number} is said to be computable if
there exists a finite, terminating algorithm such that given any positive integer $k$,
it outputs the first $k$ digits of that number\footnote{
We remark that the modern definition of computable numbers (\cite{Rice}) differs from Turing's original one.}.
 To be clear, we understand that the given object is not explicitly stated as its identification element of $\mathbb{R}$, but  means a representation or disguise that is  implicitly determined as a root of an equation (such as $\sqrt{2}$, the positive root of $x^2-2=0$), an arithmetical operation of numbers (such as $\pi+\sqrt{2}$), a Dedekind cut, a supremum,
 the limit of a convergent sequence or series (such as Euler's number $e$, see  \cite[p. 43]{Courant}), and so on.

 %, but not explicitly
%stated as an  infinite decimal (such as the Champernowne constant
%$0.12345678910\cdots$,
%obtained by concatenating the decimal representations of all natural numbers in order).  Otherwise, it is clearly computable.

%It is a surprise that Turing proved that in the sense of
%cardinality most real numbers are uncomputable.

We claim that
if it is known in advance that the sum, subtraction, product, or division of two computable numbers does not terminate,
then the arithmetical output is computable. The first three cases follow on from the  corresponding addition, subtraction and multiplication definitions. To study the division case, we can assume without loss of generality that $0<x\leq 1\leq y$.
Note for any  positive integer $k$, $(x/y_k)_k=(x_{2k}/y_k)_k$ and
\[0\leq\frac{x}{y_k}-\frac{x}{y}<10^{-k}.\]
Thus if $\theta_k({x_{2k}}/{y_k})>0$,
then $({x}/{y})_{k-1}=({x_{2k}}/{y_k})_{k-1}$. We now have two cases to consider.
In the first case, suppose there exists an increasing sequence of positive integers $\{k_i\}_{i=1}^{\infty}$
such that $\theta_{k_i}({x_{2k_i}}/{y_{k_i}})>0$ for all $i\in\mathbb{N}$.
According to the above analysis, $x/y$ is computable. In the second case, suppose there exists a non-negative integer $m$
such that $\theta_k({x_{2k}}/{y_k})=0$ for all $k>m$. Similar to the illustrating example of Definition \ref{defn12}
or the rigorous justification of Definition \ref{multiplication}, one can  show that $x/y$ is a terminating decimal.
This suffices to conclude the proof of the claim.

As an application, $\pi+\sqrt{2}$, $\pi-\sqrt{2}$, $\pi\times\sqrt{2}$ and $\pi/\sqrt{2}$  are computable because they are transcendental.
Brannan  asked a similar question about computing  $\pi+\sqrt{2}$, and his solution relies on the least upper bound property (\cite[p. 30]{Brannan}).

It is not clear which case of Hua's definition applies to the sum $\pi+e$. In fact,
it is widely believed that  $\pi+e$, $\pi-e$, $\pi e$, $\pi/e$ are all irrational (\cite{Silverman}),
and  every element of $\mathbb{Z}_{10}$  appears infinitely often in the digit sequences of $\pi$ and $e$, but a proof remains elusive.

%With Turing's concept available, it seems that our  definitions are impractical since as an example, a machine computation\footnote{
%One may assume this machine  is a sufficiently powerful calculator.} of the first digit of the sum
%\[0.4232323\cdots+0.4767676\cdots\]
%will continue to run forever.

%The authors don't know the answer to the first sum, and this is not the fault of Hua's definition.

%We are far from knowing the possible normality of $\pi$ and $e$ (\cite{Chung}),
%The answer to the second sum is the second case because this sum is an irrational number.

 \section{Conclusion}

 From Dedekind and Cantor's era  to the present day, numerous mathematicians have continuously called for
 a convincing decimal construction of the real number system.
This article, focusing on terminating decimals and elementary arithmetic rather than rational numbers and derived properties, together
with   theoretical comparison and computational discussions, is bound to have accomplished
 their wish, and provides an ideal   way for younger generations  to understand
 real numbers and their properties in the future.

\textbf{Acknowledgements}. Both authors would like to thank Prof. Jean-Pierre Demailly
for encouraging comments and helpful suggestions.
The second listed author would also like to thank Dr. Jiyou Li and Prof. Yaokun Wu for useful discussions.

\end{document}